\documentclass{article}
\usepackage{graphicx,color}
\graphicspath{{./figs/}}

\usepackage{amsthm}
\usepackage[english]{babel}

\newtheorem{theorem}{Theorem}[section]
\newtheorem{property}{Property}[section]

\usepackage{amsmath}
\usepackage{amssymb}
\usepackage{bbm}

\renewcommand{\P}[2][]{\mathbb{P}_{#1} \!\left(#2\right)}
\newcommand{\E}[2][]{\mathbb{E}_{#1} \!\left[#2\right]}

\title{Renewal Processes Represented as Doubly Stochastic Poisson Processes}
\author{Xinlong Du, Harsha Honnappa}
\date{\today}

\begin{document}

\maketitle

\begin{abstract}
This paper gives an elementary proof for the following theorem: a renewal process can be represented by a doubly-stochastic Poisson process (DSPP) if and only if the Laplace-Stieltjes transform of the inter-arrival times is of the following form:
$$\phi(\theta)=\lambda\left[\lambda+\theta+k\int_0^\infty\left(1-e^{-\theta z}\right)\,dG(z)\right]^{-1},$$
for some positive real numbers $\lambda, k$, and some distribution function $G$ with $G(\infty)=1$. The intensity process $\Lambda(t)$ of the corresponding DSPP jumps between $\lambda$ and $0$, with the time spent at $\lambda$ being independent random variables that are exponentially distributed with mean $1/k$, and the time spent at $0$ being independent random variables with distribution function $G$. 
\end{abstract}

\section{Introduction}
Doubly stochastic Poisson processes (DSPP), proposed by Cox~\cite{cox}, generalize time-inhomogeneous Poisson processes by allowing the intensity of the Poisson process to be stochastic. Thus, a DSPP can be viewed as a `mixture' point process~\cite{grandell}. The genesis of this paper is the question of what classes of counting processes can be represented as DSPPs? In other words, for a given counting process, does there exist a stochastic intensity process such that the former can be represented as a DSPP? To the best of our knowledge,~\cite{kingman} is the only paper that partially addresses this question, wherein J. F. C. Kingman proved that a stationary renewal process can be represented as a (stationary) DSPP if and only if the inter-arrival time distribution $F$ of the former has a Laplace-Stieltjes transform $\phi(\theta) = \int e^{-\theta x} dF(x)$ of the form 
\begin{align}
    \phi(\theta) = \lambda \left [ \lambda + \theta + k \int_0^\infty \left(1 - e^{-\theta z} \right)dG(z)\right]^{-1},
\end{align}
for positive real numbers $\lambda, k$ and a distribution function $G$ with $G(\infty) = 1$. Furthermore, \cite{kingman} shows that the underlying stochastic intensity of the DSPP is necessarily a Semi-Markov jump process, that jumps between the levels $\lambda > 0$ and $0$, with the time spent at level $\lambda$ being independent exponential random variables with mean $1/k$ and the time spent at $0$ being independent random variables with distribution function $G$. 

Kingman proved this remarkable theorem by appealing to the Levy-Khintchine representation \cite{kendall} for the distribution of a {\it Poisson sampled} or {\it thinned} point process $N^*$ with stationary independent increments. In particular, recall \cite{kingman} that for a process with non-negative increments, the Levy-Khintchine representation implies that there exists a non-negative number $h$ and a non-decreasing function $H$ (with support $\mathbb R_+ := [0,\infty)$) which satisfies $\int_{\mathbb R_+} z\wedge 1 dH(z) < +\infty$, such that the increments of the point process  satisfy
\begin{align}
    \mathbb E\left[ \exp\left( -\theta(N^*(x+y) - N^*(x)) \right) \right] = \exp \left( -y \left(h\theta + \int_{\mathbb R_+} (1-e^{-\theta z}) dH(z)\right) \right).
\end{align}

In this note, we provide a much simpler proof of the main theorem in~\cite{kingman} without appealing to the Levy-Khintchine theorem, using elementary arguments.

\section{Notations and Definitions}
In the sections that follow, we use $F$ to denote the distribution function of the inter-arrival times of the renewal process, and let
$$\phi_F(\theta)=\int e^{-\theta z}\,dF(z),~\theta \in [-h,\infty)~\text{for some}~ h > 0$$
be its Laplace-Stieltjes transform. Assume that $F$ is differentiable and let $f$ denote its density function\footnote{In fact, it is only necessary to assume differentiability of $F$ near $0$.}. We use $\Lambda(t)$ to denote the underlying stochastic intensity of a DSPP. 

\section{The Intensity Process}
As DSPPs are uniquely determined by the underlying intensity process, we first characterize the intensity process of a DSPP that also represents a renewal process. To achieve this, consider a stationary renewal process with the inter-arrival times distributed according to $F$.
~Fix $t<t_1<t_2<\cdots<t_n$, and consider the sequence 
$$t_{i,j}=\begin{cases}
t_i & \text{for $j=0$,}\\
\frac{1}{2}(t_{i,j-1}-t) & \text{for $j\ge 1$,}
\end{cases}$$
which converges to $t$ for all $i$ as $j\to\infty$.
Let $p_j$ denote the probability that an event occurs in each one of the intervals $[t_{i,j},t_{i,j}+dt_{i,j})$ for all $1\le i\le n$ and no other event occurs in between. For the stationary renewal process we have,
\begin{equation}\label{renewal}
p_j=\frac{1}{\mu}dt_{1,j}\prod_{i=2}^{n}f(t_{i,j}-t_{i-1,j})dt_{i,j}.
\end{equation}
On the other hand, the probability $p_j$ for the DSPP can be written as
\begin{equation}\label{dspp}
				p_j=\E{\prod_{i=1}^n\Lambda(t_{i,j})dt_{i,j}\prod_{i=1}^{n-1}e^{-\int_{t_{i,j}}^{t_{i+1,j}}\Lambda(t)\,dt}},
\end{equation}
where the expectation is taken with respect to the stochastic intensity $\Lambda$. Equate Eq.~(\ref{renewal}) and (\ref{dspp}) and take the limit as $j\to\infty$, we get
$$\frac{1}{\mu}f(0)^{n-1}=\E{\Lambda(t)^n},$$
where we have assumed right continuity of $f$ at $0$ to get the LHS, and applied dominated convergence theorem on the RHS. Now, observe that
\begin{align*}
\E{\Lambda(t)^2\left(\Lambda(t)-f(0)\right)^2}&=\E{\Lambda(t)^4-2\Lambda(t)^3f(0)+\Lambda(t)^2f(0)^2}\\
&=\frac{1}{\mu}f(0)^3-2\frac{1}{\mu}f(0)^3+\frac{1}{\mu}f(0)^3\\
&=0.
\end{align*}
Therefore, it follows that $\Lambda(t)$ equals $f(0)$ or $0$ with probability $1$. 

Next, we observe that $f(0)>0$. This can be seen from two different perspectives: consider a DSPP with the binary intensity (as described above). This process becomes `trivial' if $f(0)=0$, i.e., no event occurs in finite time, suggesting that $f(0) > 0$. We can also interpret this from the perspective of the renewal process. Note that if $F'(0)=f(0)=0$, then the probability that an event occurs within time $t$ of the previous event is (by definition of differentiation):
$$F(t)=o(t).$$
On the other hand, recall that the probability $F(t)$ that an event occurs in the interval $(0, t)$ for a DSPP satisfies
$$F(t)=\E{\Lambda(0)t}=\Theta(t)$$
when $t$ is small, suggesting that for a renewal process to be equivalent to a DSPP, $F'(0)$ must be strictly positive. This is summarized below:

\begin{property}\label{binary}
The underlying intensity process of the DSPP only takes values $\lambda>0$ and $0$. 
\end{property}

The DSPP needs to reset right after the occurrence of an event, i.e., the time before the intensity drops down to $0$ must have the same distribution regardless of when the event occurs. The only distribution that possesses this ``memoryless'' property is the exponential distribution. This is summarized as follows:
\begin{property}\label{exponential}
The amount of time the underlying intensity process of the DSPP stays at $\lambda$ is exponentially distributed with mean $1/k$, for some $k>0$.
\end{property}
The only degree of freedom left is the amount of time the intensity process stays at $0$, for which we assume some distribution function $G$ with $G(\infty)=1$. The main theorem is presented below and proven in section \ref{proof}.
\begin{theorem}\label{mainthm}
A renewal process is equivalent to a doubly-stochastic Poisson process (DSPP) if and only if the Laplace-Stieltjes transform of the inter-arrival times is of the form:
$$\phi(\theta)=\lambda\left[\lambda+\theta+k\int_0^\infty\left(1-e^{-\theta z}\right)\,dG(z)\right]^{-1},$$
for some positive real numbers $\lambda, k$, and some distribution function $G$ with $G(\infty)=1$. Furthermore, the corresponding DSPP satisfies Property \ref{binary} and \ref{exponential}.
\end{theorem}

\section{Proof of Theorem \ref{mainthm}}\label{proof}

\subsection{Necessary conditions}
For a DSPP to represent a renewal process, as argued in the previous section, it is necessary and sufficient that the intensity process $\Lambda(t)$ of the DSPP satisfies properties \ref{binary} and \ref{exponential}. Therefore, assuming any distribution function $G$ for the amount of time the stochastic intensity process $\Lambda$ stays at $0$, there exists an equivalent renewal process with some distribution function $F$. It suffices to identify the Laplace-Stieltjes transform, $\phi_F(\theta)$, of the inter-arrival times for the corresponding renewal process.

Observe that there are exactly two possibilities after the occurrence of an event: either the intensity process drops down to $0$ before another event occurs, or vice versa. Consider the former case as a ``failure'' and the latter as a ``success''. A failure occurs with probability $\frac{k}{\lambda+k}$, while a success has probability $\frac{\lambda}{\lambda+k}$; recall that $\lambda = f(0)$ and $k$ is the reciprocal of the mean time the DSPP stays at level $\lambda$. In the failure case, the intensity drops and stays at $0$ for some time before coming back to $\lambda$, after which the process ``restarts''.  

Then, if $N$ denotes the number of failures before the first success, we have
$$N\sim\text{ Geom}\left(\frac{\lambda}{\lambda+k}\right)-1,$$ where Geom$(\theta)$ denotes a geometrically distributed random variable with parameter $\theta$. At the $i$th failure, let $X_i$ be the amount of time the intensity stays at $\lambda$ before it drops to $0$ without any event occurring, and $Y_i$ the amount of time it stays at $0$ before it jumps back to $\lambda$. After $N$ rounds of failures, let $X_0$ denote the amount of time $\Lambda$ stays at $\lambda$ before an event occurs. We have
\begin{align*}
&X_i,X_0\sim^{i.i.d.} \text{ Exp}(\lambda+k)\\
&Y_i\sim^{i.i.d.} \text{ }G.
\end{align*}
Let $T$ denote the time to the next event (or the inter-arrival time), then clearly
$$T=\sum_{i=1}^{N}(X_i+Y_i)+X_0,$$ 
where the convention is taken such that the summation is $0$ if $N=0$.
The Laplace-Stieltjes transform of $T$ is\footnote{Note that we assume $\sum_{i=1}^{n-1} (X_i + Y_i) = 0$ when $n=0$.},
\begin{align*}
\E{e^{-\theta T}}=&\E{e^{-\theta X_0}}\E{e^{-\theta\sum_{i=1}^{N}(X_i+Y_i)}}\\
=&\frac{\lambda+k}{\lambda+k+\theta}\sum_{n=0}^\infty\E{e^{-\theta\sum_{i=1}^{n}(X_i+Y_i)}\Big|N=n}\P{N=n}\\ 
=&\frac{\lambda+k}{\lambda+k+\theta}\sum_{n=0}^\infty\phi_G(\theta)^n\left(\frac{\lambda+k}{\lambda+k+\theta}\right)^n\left(\frac{k}{\lambda+k}\right)^n\frac{\lambda}{\lambda+k}\\
=&\frac{\lambda}{\lambda+k+\theta-k\phi_G(\theta)}.\\
\end{align*}
We therefore conclude that the corresponding renewal process has a distribution function $F$ whose Laplace-Stieltjes transform looks like:
\begin{equation}\label{form}
\phi_F(\theta)=\E{e^{-\theta T}}=\lambda\left[\lambda+\theta+k\int_0^\infty\left(1-e^{-\theta z}\right)\,dG(z)\right]^{-1}. 
\end{equation}
In other words, for a renewal process that is identifiable as a doubly-stochastic Poisson process, it is necessarily the case that the Laplace transform of its inter-arrival time takes the form in Eq.~(\ref{form}).

\subsection{Sufficient conditions}
The sufficiency of the Laplace transform is straightforward: once we know $\phi_F(\theta)$ is of the form in Eq.~(\ref{form}), we know that the equivalent DSPP is a jump process with two states $\lambda$ and $0$. The times spent at $\lambda$ are independent exponential random variables with mean $1/k$, and the times spent at $0$ are independent random variables with distribution function $G$.  

\section{Remarks}
Assuming that both the renewal process and the corresponding DSPP are stationary, then the time before the occurrence of the first event, $t_0$, is called the ``residual time''. As a sanity check, we can verify that the residual times for the two process are also identical. From Eq.~(\ref{form}), we also get:
\begin{align*}
\phi_G(\theta)=&1+\frac{\lambda+\theta}{k}-\frac{\lambda}{k\phi_F(\theta)}\\
\phi_G'(\theta)=&\frac{1}{k}-\frac{-\lambda\phi_F'(\theta)}{k\phi_F(\theta)^2}\\
\phi_G'(0)=&\frac{1}{k}+\frac{\lambda\phi_F'(0)}{k}.
\end{align*}
For a renewal process with distribution function $F$, the Laplace transform of the residual time is known to be given by:
\begin{equation}\label{transform}
\E{e^{-\theta t_0}}=\frac{1-\phi_F(\theta)}{-\theta\phi_F'(0)}.
\end{equation}
For the corresponding DSPP,  we can divide into two cases: the process starts at $\Lambda(0)=0$ or $\Lambda(0)=\lambda$.  Using law of total expectations,
\begin{align*}
\E{e^{-\theta y_0}}=&\E{e^{-\theta y_0}|\Lambda(0)=0}\P{\Lambda(0)=0}+\E{e^{-\theta y_0}|\Lambda(0)=\lambda}\P{\Lambda(0)=\lambda}\\
=&\frac{(1-\phi_G(\theta))\phi_F(\theta)}{-\theta\phi_G'(0)}\frac{-\phi'_G(0)}{1/k-\phi'_G(0)}+\phi_F(\theta)\frac{1/k}{1/k-\phi'_G(0)}\\
=&\frac{k(1-\phi_G(\theta))\phi_F(\theta)+\theta\phi_F(\theta)}{\theta(1-k\phi'_G(0))}\\
=&\frac{1-\phi_F(\theta)}{-\theta\phi_F'(0)},
\end{align*}
which is identical to Eq.~(\ref{transform}), as desired. 

In the case when $G$ is a step function at $0$, i.e., the intensity process $\Lambda$ equals $\lambda$ always, the DSPP reduces to the regular Poisson process. It is straightforward to verify that the Laplace-Stieltjes transform of the inter-arrival times satisfies
$$\phi_F(\theta)=\lambda\left[\lambda+\theta+k\int_0^\infty\left(1-e^{-\theta z}\right)\,dG(z)\right]^{-1}=\frac{\lambda}{\lambda+\theta},$$
which is the usual exponential random variables.

\end{document}